\documentclass[twocolumn,preprintnumbers,amsmath,amssymb]{revtex4}
\usepackage[latin1]{inputenc}
\usepackage[portuges]{babel}
\usepackage[dvips]{graphicx}

\begin{document}

\title{GENERALIZATION OF THE LOGARITHM FUNCTION AND OF THE EXPONENTIAL FUNCTION WITH ARBITRARY BASE}

\author{Victor Vizcarra}

\affiliation{Departamento de Física\\Universidade Estadual de Maringá\\
Av. Colombo, 5790, Maringá - Pr, 87020-900, Brasil.}

\date{\today}

\begin{abstract}
The logarithm function and the exponential function are, by nature, base dependent. Thus, in this paper I introduces an arbitrary base in the logarithm and exponential functions, both dependent on $q$, in order to have $\log_a(x;q)$ and $a_q^x$. Some of the properties of these functions had been analyzed. The logarithm function was applied to the entropy which resulted in the $S_q = k[1 - \sum_{i = 1}^Wp_i^{\,q}]/[1 - e^{1 - q}]$ expression.
\end{abstract}

\maketitle

\section{Introduction}
The generalized entropy \cite{tsallis} is the initial point of the Non-extensive Statistical Mechanics. The base for its construction is the generalized logarithms in terms of a real parameter, $q$, inspired by concepts of multifractality. This logarithm shows a power law in its definition. In analogy with the Boltzmann-Gibbs entropy, the general entropy must contain the generalized neperian logarithm. It happens that the definition of this function did not find any remnants of any particular base; obviously, is the neperian base, but this base is not explicit. Why? The same problem, if it can be called a problem, is carried by the exponential function. Both functions, the logarithm function and the exponential function, are by nature base dependent. In this work I propose to define the functions logarithm and exponential in order to load explicitly the base on which they evolve, the aim is just to do a mathematical analysis.

\section{The Logarithm Function}
As mentioned before, the logarithm function should carry an explicit basis. This set me to function as

\begin{equation}
\log_a(x;q) = \displaystyle\frac{x^{1 - q} - 1}{a^{1 - q} - 1}
\label{log}
\end{equation}
where $a$ represents the base of the logarithm. The traditional logarithm is recovered when $q = 1$.
To avoid confusion, the generalized logarithm, or $q$-logarithm, will be shown only as a logarithm; when the case is $ q = 1$, the logarithm will be treated as traditional logarithm or usual logarithm. The same attitude will be taken when dealing with the exponential function.

The Eq. (\ref{log}) say categorically that the logarithm follows a law power when $q \neq 1$. For any value of $q$, easily be seen the essential characteristics of the traditional logarithm, to know, the logarithn of one is always zero, i.e., $\log_a(1;q) = 0$, $\forall q,a \neq 1$; the logarithm of the base, whatever the basis, is always one: $\log_a(a;q) = 1$, $\forall q,a \neq 1$.

Another characteristic easily verifiable is the change of ownership of base
\begin{equation}
\log_a(x;q) = \displaystyle\frac{\log_b(x;q)}{\log_b(a;q)}.
\end{equation}

Certainly, we must expect winning and losing properties simply because we changed the definition. It is very stimulating do touch on properties that exist only in areas prohibited in the usual logarithm, see for example, the logarithm of a negative number, in this case, the value is real and finite for any $q$ integer, different from 1.
\begin{equation}
\log_a(-x;q) = \displaystyle\frac{1 + (-1)^q\,x^{1 - q}}{1 - a^{1 - q}};\ \ \ \ \ \ (q \in \mathbb{Z} - \{1\}).
\end{equation}
The logarithm of zero also presents real and finite result to certain $q$ conditions.
\begin{equation}
\log_a(0;q) = \displaystyle\frac1{1 - a^{1 - q}} ;\ \ \ \ \ \ (q < 1),
\end{equation}
In another case, $(q \geq 1)$, the value is infinite, $\pm\infty$. Here, the positive or negative sign depends on the value of $a$.

Another interesting result is the base logarithm equals zero, $a = 0$. Its value is real and finite for $q$ special values, in the region $(q < 1)$,
\begin{equation}
\log_0(x;q) = 1 - x^{1 - q} ;\ \ \ \ \ \ (q < 1),
\end{equation}
in the region$(q \geq 1)$, $\log_0(x;q)$ if void.

In the traditional version of the logarithm are considered absolute three properties that, unfortunately, in the new definition are modified, if not lost. The power law, for example, is violated, $\log_a(x^k;q) \neq k\cdot\log_a(x;q)$, but, a negative power can be rearranged in order to obtain
\begin{equation}
\log_a(x^{-r};q) = -\displaystyle\frac1{x^{r\cdot(1 - q)}}\cdot\log_a(x^r;q).
\end{equation}
The multiplation law also changes, it results in a pseudo-additivity which not only depends on the value of $q$, but also the base $a$, i.e.,
\begin{eqnarray}
\!\!\!\!\!\log_a(x\cdot y;q) &=& \log_a(x;q) + \log_a(y;q)\nonumber\\[0.1cm]
&& +\ (a^{1 - q} - 1)\log_a(x;q)\log_a(y;q),
\end{eqnarray}
and similarly, the ratio of two variables,
\begin{eqnarray}
\!\!\!\!\!\!\!\!\log_a\left(\displaystyle\frac{x}{y};q\right) &=& \log_a(x;q) - \displaystyle\frac1{y^{1 - q}}\log_a(y;q)\nonumber\\[0.1cm]
&& -\ \displaystyle\frac{(a^{1 - q} - 1)}{y^{1 - q}}\log_a(x;q)\cdot\log_a(y;q).
\end{eqnarray}

The assymptotic behavior of the logarithm (\ref{log}) also presents results that differ from the usual logarithm for certain values of $q$. For example, in the $q > 1$ region, within the limits $x \rightarrow \infty$ and $a \rightarrow \infty$, the logarithm provides us with finite values; the usual result is recovering in the region $q \leq 1$, see below,
\begin{eqnarray}
&\displaystyle\lim_{x\rightarrow\infty}\log_a(x;q) = \left\{\begin{array}{ll}
\displaystyle\frac1{1 - a^{1 - q}}\ ;& (q > 1),\\[.3cm]
\ \ \ \ \infty\ ;& (q \leq 1),
\end{array}\right.\\[.3cm]
&\displaystyle\lim_{a\rightarrow\infty}\log_a(x;q) = \left\{\begin{array}{ll}
1 - x^{1 - q}\ ;& (q > 1),\\[0.3cm]
\ \ \ \ 0\ ;& (q \leq 1).
\end{array}\right.
\end{eqnarray}

\begin{figure}[htb]
\includegraphics[width=10.cm,height=8.cm]{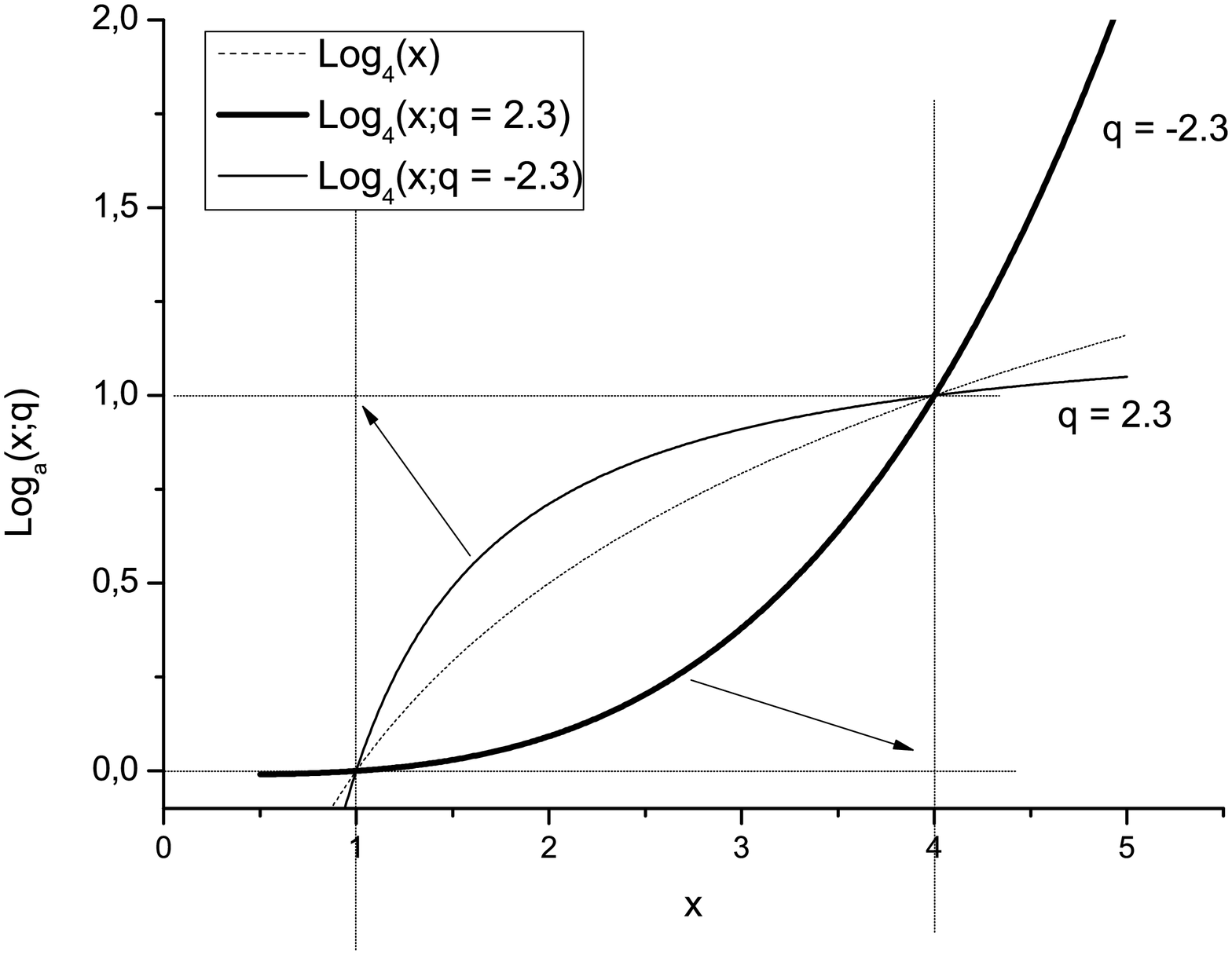}
\caption{\label{fig.1}Graphic representation of the logarithm, for fixed values of $q = 2.3$ and $q = -2.3$, and $a = 4$. The arrows indicate the assymtotic behavior of the function when $q \rightarrow \infty$.}
\end{figure}

\begin{figure}[h]
\includegraphics[width=10.cm,height=8.cm]{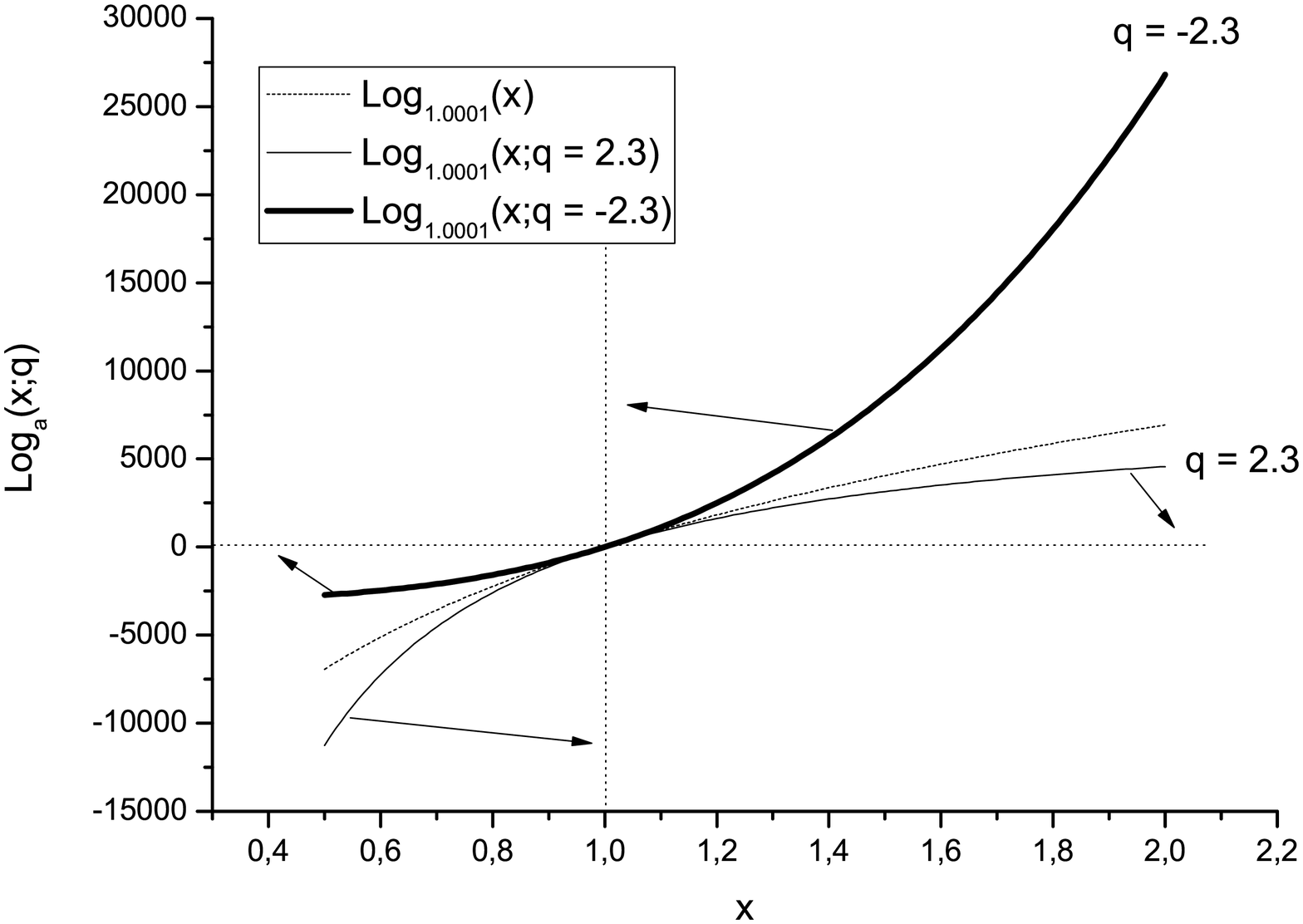}
\caption{\label{fig.2}Graphic representation of the logarithm, for fixed values of $q = 2.3$ and $q = -2.3$, and $a = 1.0001$. The arrows indicate the assymtotic behavior of the function when $q \rightarrow \infty$.}
\end{figure}

In the Figures \ref{fig.1} and \ref{fig.2} are represented the generalized logarithm graphics in $q$ and $a$, when $a$ is the base, fixed in two particular values for these parameters: $q = 2.3$ and $q = -2.3$, and $a = 1.0001$ and $a = 4$. The traditional logarithm is also plotted in both figures, they correspond to the dashed line. In the figure \ref{fig.1}, the $q$-logarithms present opposing concavities; the distance between these concavities increases as the amount of $q$ also increases, but, tending to the limits $\log_a(x;-\infty) = 1$ and $x = a$ when $q \rightarrow \infty$.

In the figure \ref{fig.2}, the concavities are no longer separated, the curves tend to touch on $x = 1$ when $a \rightarrow 1$. Even under these conditions, the curvature of the lines will also be close to its maximum when $q \rightarrow \infty$, and these limits are $\log_a(x;-\infty) = 1$ and $x = a$, in the case, $a = 1.0001$.

It is worth to mention that the tendency to such limits is a feature itself of logarithms. When $q \rightarrow -\infty$ the logarithm grows to infinity but tending to value of the base, i.e., $x = a^+$ (the superscript means the tending from the right-hand side toward to the left-hand side). For $q \rightarrow +\infty$, the values of the logarithm decrease negatively near $x = 1$. All this happens for arbitrary values of base more than 1, $a > 1$. This same feature is also where the base takes values between 0 and 1, $0 < a < 1$, but, on the left side of $x = 1$.

At this point we must call the attention to the notation used by all for the natural logarithm and the 10 base logarithm. In agreement with the universal notation, we wrote $\log(x;q)$ for the 10 base logarithm and $\ln(x;q)$ for natural logarithm, corresponding to the $e = 2.718281828...$ base. The logarithm proposed by Tsallis \cite{tsallis} is derived from the natural logarithm, being enough only to expand the base factor $e^ {1 - q} - 1$ in $1 - q$ potencies, i.e.,
\[e^{1 - q} - 1 = (1 - q) + \frac{(1 - q)^2}{2!} + \frac{(1 - q)^3}{3!} + \frac{(1 - q)^4}{4!} + \cdots\]
Since the series presents infinite terms, it is obvious that we have infinite expressions approximated for the generalized logarithm. For example, leaving only the first approximation of expansion, which is the linear term, we have exactly the logarithm of Tsallis. Another way to approximate logarithm would replace the base factor of the quadratic approximation of the expansion, i.e., $(1 - q)(3 - q)/2$. Thus, we can get infinite expressions, an expression for each truncated term of the expansion.

Of course, we cannot simply disdain the nonlinear terms. Strictly, whatever the base us to use, the correct term for the logarithm that is defined in Eq. (\ref{log}), since it involves not just the linear term in $1 - q$, but every non-linear terms of the expansion $e^{1 - q} - 1$.

\section{The Exponential Function}
The exponential function also expresses the fact of being dependent on base. The defined as
\begin{equation}
a^x_q = [1 + (a^{1 - q} - 1)\cdot x]^{1/(1 - q)}
\label{exp}
\end{equation}
where $a$ represents the base of the exponential. We can verify two essential characteristics easily being preserved for this definition, exponential of the unit and exponential of zero: $a^1_q = a,\ \forall q$, and $a^0_q = 1,\ \forall a,q$.

A curious consequência proceeding from the Eq. (\ref{exp}) is the nonzero value resulting from the exponential of zero base, result that only happens in the $q < 1$ region:
\begin{equation}
0^x_q = (1 - x)^{1/(1 - q)} ;\ \ \ \ \ \ (q < 1)
\end{equation}
outside that region, i.e., $q \geq 1$, the exponential disappears.

As well as in the case of the logarithm, the exponential function also violates some of its properties. For example, the exponential of the sum of variables is not the product of the exponential, i.e., $a^{x + y}_q \neq a^x_q\cdot a^y_q$. The term of the left-hand side of the inaquality is written as
\begin{equation}
a^{x + y}_q = [(a^x_q)^{1 - q} + (a^y_q)^{1 - q} - 1]^{1/(1 - q)}
\label{exp_da_soma}
\end{equation}
and the term of the right-hand side, written as
\begin{equation}
a^x_q\cdot a^y_q = a_q^{x + y + (a^{1 - q} - 1)xy}.
\end{equation}
The ratio of two exponential is not equal to the difference of exponential, first because the expression of ratio is
\begin{equation}
\displaystyle\frac{a^x_q}{a^y_q} = a_q^{(x - y)/[1 + (a^{1 - q} - 1)y]}
\end{equation}
and the difference is decurrent of the Eq. (\ref{exp_da_soma}). Another property which is not fulfilled is the power of an exponential, i.e., $[a_q^x]^n \neq a_q^{nx}$.

\begin{figure}[htb]
\includegraphics[width=10.cm,height=8.cm]{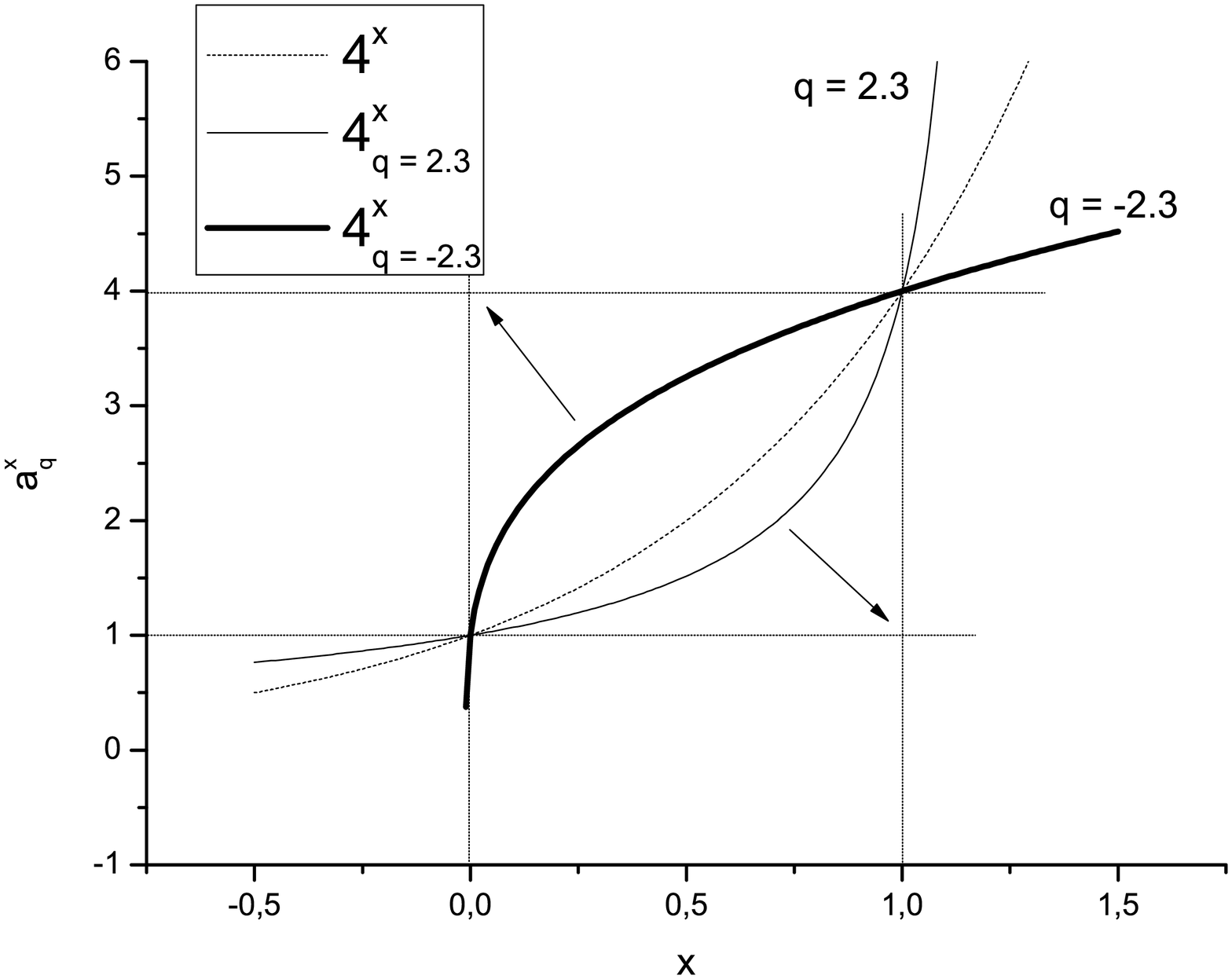}
\caption{\label{fig.3}Graphic representation of the logarithm, for fixed values of $q = 2.3$ e $q = -2.3$, and $a = 4$. The arrows indicate the assymtotic behavior of the function when $q \rightarrow \infty$.}
\end{figure}

In the figure \ref{fig.3} we plotted the graph of the exponential function for two fixed values of $q$: $q = 2.3$ and $q = -2.3$, with the base value equal to the $a = 4$.

The dashed line is the exponential traditional and the other two, corresponding to different values of $q$. Since the value of the base is a positive number, $a > 0$, and the different of the unit, the concavities of the curves will be always opposing, and the greater the amount of $q$, in such a way more the bendings will be more accented, the limit is the value of their base in ordinate axis, this in the case of $q < 0$, and for the case $q > 0$, the limit is 1 in abscissa axis.The exponential, however, obeys the restriction $x > 1/(1 - a^{1 - q})$ for $q > 0$ and $x > 1/(1 - a^{1 + q})$ for $q < 0$.

When the value of the base is close to the unit, $a \rightarrow 1$, the curves tend to touch on $x = 1$, following the same rules addressed in the preceding paragraph.

\section{APPLICATION TO ENTROPY}
In this paper does not intend to do some study on the generalized entropy \cite{tsallis}; I just want to show that, despite the contributions of accurate Nonextensive Statistical Mechanics in favor of the experimental results, the correct expression of entropy is not the one that has been used until the moment, and yes the decurrent one of the logarithms of the expression (\ref{log}). A entropia de Boltzmann-Gibbs está assentada sobre o logaritmo neperiano, de base igual a $a = e = 2.718281828\ldots $. In analogy to the entropy of Boltzmann-Gibbs, $S = k\sum_{i = 1}^Wp_i\ln(1/p_i)$, we wrote the generalized entropy replacing the neperian logarithm by logarithm of the expression (\ref{log}) with equal base $a = e = 2.718281828\ldots $, i.e., 
\[S_q = k\sum_{i = 1}^Wp_i\ln\left(\dfrac1{p_i};q\right)\]
which, of course, will result in
\begin{equation}
S_q = k\,\displaystyle\frac{1 - \sum_{i = 1}^Wp^q_i}{1 - e^{1 - q}}\,,\ \ \ \ \ \ (q \in \mathbb{R})
\label{entropia_victor}
\end{equation}
bounded to the condition $\sum_{i = 1}^Wp_i = 1$, where $W \in \mathbb{N}$ represents the total number of possible microscopical configurations of a system in special and $\{p_i\}$ represents the probabilities associates. We may notice that the entropy proposed by Tsallis \cite{tsallis} just corresponds to the linear term in $1 - q$ of the $e^{1 - q} - 1$ expansion. Then, the above equation is the correct form of the generalized entropy, since it involves not just the linear term in $1 - q$, but every non-linear terms of the expansion $e^{1 - q} - 1$.

An entropy with the equal base $a = e$, like Eq. (\ref{entropia_victor}), it is not absolute, because the base can be any value, depending on the system that one is studying. For example, another particular expression of the entropy was obtained by Daróczy \cite{daroczy,aczel}, studying measures parameters of the information. The base parameter of the entropy of Daróczy is $a = 2$, a value close to $e$-neperian.

I mean that other entropy can be obtained using the logarithm of expression (\ref{log}), just simply manipulate, mutilate and eliminate any terms of his own expression, and like this to build a statistical mechanics for that entropy.

\section{Conclusion}
In this paprer I have defined the logarithm function and the exponencial exponential in order to respect the nature primary of these functions, which is the total dependence on base. The objective of this work treated only the mathematical side of these functions, since mathematical 'coherent' functions they are useful tools for any study of a natural system. The initial application of the generalized incomplete logarithm, and what it gave rise to the whole field of the Nonextensive Statistical Mechanics, was the entropy \cite{tsallis}. Now, applying the full logaruthm of Eq (\ref{log}) to the entropy, certainly we will be able to understand in more accurate those divested situations of the stationary states when were known as a thermal equilibrium. It is clear that the functions (\ref{log}) and (\ref{exp}) is not restricted only to the field of thermodynamics but also to other areas in which any state of non-equilibrium is common, as you can see, for example, in economic sciences.

Finally, I must mention which the construction of a physical theory requires of us the ability to we can write any kind of equations that can describe any phenomenon that we can study, worth is, of course, of its mathematical properties and of his logical handlings. So, the theory what we build must emerge solidly and must grow that it releases of any type of ambiguity. For such a reason, anyone function what we claim to use and to generalize, it must guarantee in absolute way that his basic properties should keep on lasting a long time while the theory is constructed. For such a reason I find unquestionable the truthfulness of the equations (\ref{log}) and (\ref{exp}).

\end{document}